\documentclass[12pt,twoside]{amsart}
\usepackage{amssymb}
\usepackage{amscd}

\title[toric Fano 5-folds]
{Smooth toric Fano five-folds of index two} 
\author{Hiroshi Sato} 
\thanks{The author is partly supported by the 
Grant-in-Aid for JSPS Fellows, The Ministry of 
Education, Science, Sports and Culture, Japan.}
\address{
Osaka City University Advanced Mathematical Institute, 
3-3-138 Sugimoto, Sumiyoshi-ku, Osaka 558-8585, 
Japan}
\email{hirosato@sci.osaka-cu.ac.jp}

\newcommand{\Hom}[0]{{\operatorname{Hom}}}

\newcommand{\G}{\mathop{\rm G}\nolimits}
\newcommand{\PC}{\mathop{\rm PC}\nolimits}
\newtheorem{thm}{Theorem}[section]
\newtheorem{lem}[thm]{Lemma}
\newtheorem{cor}[thm]{Corollary}
\newtheorem{prop}[thm]{Proposition}
\newtheorem{lem-def}[thm]{Definition-Lemma}

\theoremstyle{definition}
\newtheorem{ex}[thm]{Example}
\newtheorem{defn}[thm]{Definition}
\newtheorem{rem}[thm]{Remark}
\newtheorem*{ack}{Acknowledgments}       
 
\begin{document}
\bibliographystyle{amsalpha+}

\begin{abstract}
In this paper, we classify smooth toric Fano 5-folds of 
index 2. There exist exactly $10$ smooth toric Fano 5-folds 
of index 2 up to isomorphisms.
\end{abstract}

\maketitle
\tableofcontents

\section{Introduction}\label{intro}

\thispagestyle{empty}

For a smooth Fano $d$-fold $X$, 
the {\em index} $i_X$ of $X$ is defined 
as follows:
$$i_X:=\max\left\{m\in{\mathbb Z}_{\geq 1}\,\left|\,
-K_X=mH\mbox{ for a Cartier divisor }H\right.\right\}.$$
There is a famous result of 
\cite{kobayashiochiai} which say $1\leq i_X\leq d+1$ and 
a smooth Fano $d$-fold of index $d+1$ or $d$ is isomorphic to 
${\mathbb P}^d$ or $Q^d$, respectively, where $Q^d$ is 
the $d$-dimensional quadric. A smooth Fano $d$-fold of index $d-1$ 
or $d-2$ is called a {\em del Pezzo} manifold or a {\em Mukai} manifold, 
respectively, and there are classifications for these 
manifolds (see \cite{fujita}, \cite{mella} and \cite{mukai}). 

So, the next problem is the classification of 
smooth Fano $d$-folds of index $d-3$. 
If $d\geq 6$ and the Picard number is greater than $1$, 
there is the classification (see \cite{wisn}). 
For the case $d=5$, there are some partial classifications 
(see \cite{novelli2} and \cite{novelli1}).
Toward the general classification, in this paper, 
we classify smooth toric Fano $5$-folds of index $2$. 
We show that there exist exactly $10$ smooth toric Fano $5$-folds of 
index $2$ (see Theorem \ref{mainresult}). 
We remark that 
since a smooth complete toric $d$-fold of Picard number $1$ is 
isomorphic to ${\mathbb P}^d$, this result completes the classification 
of smooth toric Fano $d$-folds of index $d-3$. 

The content of this paper is as follows: Section \ref{junbi} is 
a section for preparation. We review the Mori theory 
for smooth toric varieties. In Section \ref{maindabesa}, 
we consider the classification of smooth toric Fano 
$5$-folds of index $2$. Section \ref{nonsplitting} is devoted to 
constructing an example of a smooth toric Fano $7$-fold of 
index $2$ which has no projective space bundle structure. 


\begin{ack} 
The author would like to thank Doctor Ichitaka Suzuno for 
advice and encouragement. He is greatful to Professor 
Masa-Nori Ishida who gave him useful comments for the 
construction of Example \ref{hanachangum}. 
\end{ack}

\section{Preliminaries}\label{junbi}

In this section, we explain some basic 
facts of the toric geometry. See \cite{bat}, 
\cite{bat4}, \cite{fujisato}, \cite{fulton}, \cite{oda} and 
\cite{sato1} for the detail. 

Let $\Sigma$ be a nonsingular complete fan in $N:=\mathbb{Z}^d$, 
$M:={\rm Hom}_{\mathbb{Z}}(N,\mathbb{Z})$ 
and $X=X_{\Sigma}$ the associated 
smooth complete toric $d$-fold 
over an algebraically closed field $k$. 
Let $\G(\Sigma)$ be the set of primitive generators of 
$1$-dimensional cones in $\Sigma$. A subset $P\subset 
\G(\Sigma)$ is called a {\em primitive collection} 
if $P$ does not generate a cone in $\Sigma$, while 
any proper subset of $P$ generates a cone in $\Sigma$. 
We denote by $\PC(\Sigma)$ the set of primitive 
collections of $\Sigma$. 
For a primitive collection $P=\{x_1,\ldots,x_m\}$, there exists 
the unique cone $\sigma(P)$ in $\Sigma$ such that $x_1+\cdots+x_m$ 
is contained in its relative interior 
since $\Sigma$ is complete. So, we obtain an equality
\begin{equation}\label{holyland}
x_1+\cdots+x_m=b_1y_1+\cdots+b_ny_n,
\end{equation}
where $y_1,\ldots,y_n$ are the generators of $\sigma(P)$, 
that is, $\sigma(P)\cap\G(\Sigma)=\{y_1,\ldots,y_n\}$, 
and $b_1,\ldots,b_n$ are positive integers. 
We call this equality the {\em primitive relation} of $P$. 
By the standard exact sequence 
$$0\to M\to\mathbb{Z}^{\G(\Sigma)}\to{\rm Pic}(X)\to0$$
for a smooth toric variety, we have 
\begin{eqnarray}
A_1(X) & \simeq & \Hom_{\mathbb Z}({\rm Pic}(X),{\mathbb Z})
\simeq\Hom_{\mathbb Z}({\mathbb Z}^{\G(\Sigma)}/M,{\mathbb Z}) 
\nonumber \\
 & \simeq & M^{\perp}\subset
\Hom_{\mathbb Z}({\mathbb Z}^{\G(\Sigma)},{\mathbb Z}),\nonumber
\end{eqnarray}
where $A_{1}(X)$ 
is the group of $1$-cycles on $S$ modulo rational 
equivalences, 
and hence
$$A_{1}(X)\simeq \left\{(b_{x})_{x\in\G(\Sigma)}\in
\Hom_{\mathbb Z}({\mathbb Z}^{\G(\Sigma)},{\mathbb Z})
\;\left| \;\sum_{x\in \G(\Sigma)} b_{x}x=0\right.\right\}.$$
Thus, by the equality $x_1+\cdots+x_m-(b_1y_1+\cdots+b_ny_n)=0$, 
we obtain an element $r(P)$ in $A_{1}(X)$ for each 
primitive collection $P\in\PC(\Sigma)$. 
We define the {\em degree} of $P$ as 
$\deg P:=\left( -K_X\cdot r(P)\right)=m-(b_1+\cdots+b_n)$. 

\begin{prop}[\cite{bat} and \cite{reid}]\label{toriccone}

Let $X=X_{\Sigma}$ be a smooth projective toric variety. Then, 
the Mori cone of $X$ is described as
$${\rm NE}(X)=\sum_{P\in\PC(\Sigma)}{\mathbb R}_{\geq 0}r(P)
\subset A_{1}(X)\otimes\mathbb{R}.$$
\end{prop}
A primitive collection $P$ is said to be 
{\em extremal} if $r(P)$ is contained 
in an extremal ray of ${\rm NE}(X)$. 

\begin{rem}\label{keitai1}
If $x_1+\cdots+x_m=b_1y_1+\cdots+b_ny_n$ is 
an extremal primitive relation, then $m+n\le d+1$, 
because $r(P)$ corresponds to an irreducible torus invariant curve.
\end{rem}

\begin{cor}\label{keitai2}
Let $X=X_{\Sigma}$ be a smooth projective toric variety. Then, 
$X$ is {\em Fano} if and only if $\deg P>0$ for any extremal primitive 
collection $P\in\PC(\Sigma)$.
\end{cor}

For extremal primitive relations, we need the following propositon 
and definition for the classification.

\begin{prop}[\cite{cas} and \cite{sato1}]\label{contr}
Let $X=X_{\Sigma}$ be a smooth projective toric variety and 
$P$ an extremal primitive collection. Then, 
for any $P'\in\PC(\Sigma)\setminus\{P\}$ 
such that $P\cap P'\neq\emptyset$, 
$$\left(P\setminus P'\right)\cup\left(\sigma(P)\cap\G(\Sigma)\right)$$
contains a primitive collection.
\end{prop}

\begin{defn}[\cite{bat}]
Let $X=X_{\Sigma}$ be a smooth complete toric variety. 
Then, $\Sigma$ is a {\em splitting} fan if 
$P\cap P'=\emptyset$ for any $P,\ P'\in\PC(\Sigma)$ such that 
$P\neq P'$.
\end{defn}

If $\Sigma$ is a splitting fan, then there exists a sequence of 
smooth complete toric varieties
$$X=X_{\Sigma}=:X_s\stackrel{\varphi_s}{\to}
X_{s-1}\stackrel{\varphi_{s-1}}{\to}\cdots
\stackrel{\varphi_3}{\to}X_2\stackrel{\varphi_2}{\to}
X_1\simeq{\mathbb P}^l,$$
where $X_i\stackrel{\varphi_i}{\to}X_{i-1}$ is a toric 
projective space bundle and $l\in{\mathbb Z}_{\geq 1}$. 
We remark that $s$ is the Picard number of $X$. 
The number of the primitive collections of $\Sigma$ is 
also $s$. 

\section{Classification}\label{maindabesa}

We start the classification. 

Let $X=X_{\Sigma}$ be a smooth toric fano $5$-fold of index $2$. 
In this case, $\deg P$ is an even number for any $P\in\PC(\Sigma)$. 
Then, by Remark \ref{keitai1} and Corollary \ref{keitai2}, 
the type of any extremal primitive relation is one of 
the following:
\begin{enumerate}
\item $x_1+x_2+x_3+x_4+x_5+x_6=0,$
\item $x_1+x_2+x_3+x_4+x_5=y_1,$
\item $x_1+x_2+x_3+x_4+x_5=3y_1,$
\item $x_1+x_2+x_3+x_4=0,$
\item $x_1+x_2+x_3+x_4=2y_1,$
\item $x_1+x_2+x_3+x_4=y_1+y_2,$
\item $x_1+x_2+x_3=y_1$ and 
\item $x_1+x_2=0,$
\end{enumerate}
where $\{x_1,x_2,x_3,x_4,x_5,x_6,y_1,y_2\}\subset\G(\Sigma)$. 

First of all, the existence of 
an extremal primitive relation of 
type $(1)$ imply that $X\cong{\mathbb P}^5$, but 
${\mathbb P}^5$ is of index $6$. 
So, there does not exist 
an extremal primitive relation of type $(1)$.

\begin{prop}\label{shinbun1}
Let $X=X_{\Sigma}$ be a smooth toric Fano $5$-fold of index $2$. 
If $X$ has an extremal primitive relation of type 
$(2)$ or $(3)$, 
then $X$ is isomorphic to either 
$${\mathbb P}_{{\mathbb P}^4}
\left(\mathcal{O}_{{\mathbb P}^4}\oplus
\mathcal{O}_{{\mathbb P}^4}(1)\right)
\mbox{ or }
{\mathbb P}_{{\mathbb P}^4}
\left(\mathcal{O}_{{\mathbb P}^4}\oplus
\mathcal{O}_{{\mathbb P}^4}(3)\right).$$
\end{prop}

\begin{proof}
In this case, $X$ has a divisorial contraction whose image 
of the exceptional divisor 
is a point. So, only we have to do is to check 
the classified list in \cite{bonavero}. 
\end{proof}

\begin{prop}\label{shinbun2}
Let $X=X_{\Sigma}$ be a smooth toric Fano $5$-fold of index $2$. 
If $X$ has an extremal primitive relation of type $(4)$, 
then $X$ is isomorphic to either 
$${\mathbb P}^1\times{\mathbb P}^1\times{\mathbb P}^3
\mbox{ or }
{\mathbb P}_{{\mathbb P}^2}
\left(\mathcal{O}_{{\mathbb P}^2}\oplus
\mathcal{O}_{{\mathbb P}^2}\oplus
\mathcal{O}_{{\mathbb P}^2}\oplus
\mathcal{O}_{{\mathbb P}^2}(1)\right).$$
\end{prop}

\begin{proof}
In this case, $X$ is a ${\mathbb P}^3$-bundle over 
a toric del Pezzo surface. By checking the classification 
of toric del Pezzo surfaces, we can prove this proposition. 
\end{proof}

\begin{prop}\label{shinbun3}
Let $X=X_{\Sigma}$ be a smooth toric Fano $5$-fold of index $2$. 
If $X$ has an extremal primitive relation of type $(5)$, 
then $X$ is isomorphic to 
$${\mathbb P}^1\times{\mathbb P}_{{\mathbb P}^3}
\left(\mathcal{O}_{{\mathbb P}^3}\oplus
\mathcal{O}_{{\mathbb P}^3}(2)\right).$$
\end{prop}

\begin{proof}
In this case, $X$ has a divisorial contraction whose image 
of the exceptional divisor 
is a curve. So, only we have to do is to check 
the classified list in \cite{sato2}. 
\end{proof}

\begin{prop}\label{shinbun4}
Let $X=X_{\Sigma}$ be a smooth toric Fano $5$-fold of index $2$. 
If $X$ has an extremal primitive relation of type $(8)$, 
then $X$ has a ${\mathbb P}^1$-bundle structure. 
In this case, there exist exactly $9$ 
such smooth toric Fano $5$-folds of index $2$ 
$($see {\rm Theorem \ref{mainresult}}$)$.
\end{prop}

\begin{proof}
By the classified list of smooth toric Fano $4$-folds 
(see \cite{bat4} and \cite{satopubl}), 
we have exactly 
$8$ smooth toric Fano $4$-folds whose indices are 
at least $2$. 
It is an easy exercise to construct toric 
${\mathbb P}^1$-bundles which are smooth toric Fano 
$5$-folds of index $2$ over them.
\end{proof}

Thus, we may assume that every extremal primitive 
relation of $X$ is of type $(6)$ or $(7)$. 
However, there is no such variety as follows.

\begin{lem}\label{shinbun5}
Let $X=X_{\Sigma}$ be a smooth toric Fano $5$-fold of index $2$. 
Then, $X$ has a primitive relation of type other than 
$(6)$ and $(7)$.
\end{lem}

\begin{proof}
Suppose that every extremal primitive 
relation of $X$ is of type $(6)$ or $(7)$. 

First of all, we claim that there is no primitive 
collection $P\in\PC(\Sigma)$ such that $\# P=2$. 
This is obvious because $P$ has to be an extremal 
primitive collection. Namely, its primitive relation is 
of type $(8)$.  

Suppose that there exists an extremal primitive relation 
$x_1+x_2+x_3+x_4=y_1+y_2$. If the Picard number of $X$ 
is two, then $X$ has at least one Fano contraction. 
So, we may assume that there exist two distinct elements 
$z_1,\ z_2\in\G(\Sigma)\setminus\{x_1,x_2,x_3,x_4,y_1,y_2\}$. 
Thus, we have two extremal primitive relations
$$y_1+y_2+z_1=w_1\mbox{ and }y_1+y_2+z_2=w_2,$$
where $w_1,\ w_2\in\G(\Sigma)$. However, Proposition \ref{contr} 
says that $\{z_1,w_2\}$ and $\{z_2,w_1\}$ are primitive collections, 
and this is a contradiction. 

Finally, we may assume that every extremal primitive 
relation of $X$ is of type $(7)$. 

As above, 
for any distinct extremal primitive collections 
$P_1,\ P_2\in\PC(\Sigma)$, we have $\#(P_1\cap P_2)\neq 2$. 
So, let $\#(P_1\cap P_2)=1$, and let $x_1+x_2+x_3=y_1$ 
and $x_1+x_4+x_5=y_2$ be the corresponding primitive relations. 
Then, Propositon \ref{contr} imply that $\{x_2,x_3,y_2\}$ is 
a primitive collection. This primitive collection is extremal. 
So, the corresponding primitive relation is $x_2+x_3+y_2=z$ 
for some $z\in\G(\Sigma)$. By applying Proposition \ref{contr} again, 
$\{x_1,z\}$ is a primitive collection. This is a contradiction. 
Therefore, $P_1\cap P_2=\emptyset$. 

Since the Picard number of $X$ is at least $3$, 
there exist at least three extremal primitive collections 
$P_1$, $P_2$ and $P_3$. Thus, $\#\G(\Sigma)\geq 9$ and 
the Picard number of $X$ is at least $4$. So, 
we have a new extremal primitive collection $P_4$ and 
$\#\G(\Sigma)\geq 12$. We can continue this process endlessly. 
This is impossible.
\end{proof}

By Propositions \ref{shinbun1}, \ref{shinbun2}, 
\ref{shinbun3}, \ref{shinbun4} and Lemma \ref{shinbun5}, 
we complete the classification:

\begin{thm}\label{mainresult}
Let $X=X_{\Sigma}$ be a smooth toric Fano $5$-folds of index $2$. 
Then, $X$ is one of the follwing$:$
\begin{enumerate}
\item ${\mathbb P}_{{\mathbb P}^2}
\left(\mathcal{O}_{{\mathbb P}^2}\oplus
\mathcal{O}_{{\mathbb P}^2}\oplus
\mathcal{O}_{{\mathbb P}^2}\oplus
\mathcal{O}_{{\mathbb P}^2}(1)\right)$.
\item ${\mathbb P}_{{\mathbb P}^4}
\left(\mathcal{O}_{{\mathbb P}^4}\oplus
\mathcal{O}_{{\mathbb P}^4}(1)\right)$.
\item ${\mathbb P}_{{\mathbb P}^4}
\left(\mathcal{O}_{{\mathbb P}^4}\oplus
\mathcal{O}_{{\mathbb P}^4}(3)\right)$.
\item ${\mathbb P}^1\times{\mathbb P}^1\times{\mathbb P}^3$.
\item ${\mathbb P}^1\times{\mathbb P}_{{\mathbb P}^3}
\left(\mathcal{O}_{{\mathbb P}^3}\oplus
\mathcal{O}_{{\mathbb P}^3}(2)\right)$.

\item  ${\mathbb P}^1$-bundle over 
${\mathbb P}_{{\mathbb P}^2}
\left(\mathcal{O}_{{\mathbb P}^2}\oplus
\mathcal{O}_{{\mathbb P}^2}\oplus
\mathcal{O}_{{\mathbb P}^2}(2)\right)$ 
whose primitive relations are
$x_1+x_2+x_3=x_4$, $x_4+x_5+x_6=x_7$ and $x_7+x_8=0$, 
where $\G(\Sigma)=\{x_1,\ldots,x_8\}$.

\item  ${\mathbb P}^1$-bundle over 
${\mathbb P}^2\times{\mathbb P}^2$ 
whose primitive relations are
$x_1+x_2+x_3=x_7$, $x_4+x_5+x_6=x_7$ and $x_7+x_8=0$, 
where $\G(\Sigma)=\{x_1,\ldots,x_8\}$.

\item  ${\mathbb P}^1$-bundle over 
${\mathbb P}^2\times{\mathbb P}^2$ 
whose primitive relations are
$x_1+x_2+x_3=x_7$, $x_4+x_5+x_6=x_8$ and $x_7+x_8=0$, 
where $\G(\Sigma)=\{x_1,\ldots,x_8\}$.

\item ${\mathbb P}^1\times{\mathbb P}^1\times
{\mathbb P}_{{\mathbb P}^2}\left(\mathcal{O}_{{\mathbb P}^2}\oplus
\mathcal{O}_{{\mathbb P}^2}(1)\right)$.
\item ${\mathbb P}^1\times{\mathbb P}^1\times{\mathbb P}^1\times
{\mathbb P}^1\times{\mathbb P}^1$.
\end{enumerate}
\end{thm}

\begin{rem}
As in Theorem \ref{mainresult}, the fan of 
every smooth toric Fano $5$-fold of 
index $2$ is a splitting fan. 
Moreover, if $d\leq 5$ and $p\geq 2$, then 
the fan of 
every smooth toric Fano $d$-fold of 
index $p$ is a splitting fan. 
In section \ref{nonsplitting}, we show 
higher dimensional examples of toric Fano 
manifolds of higher indicies which admit no 
projective space bundle structure. 
\end{rem}

\section{Example}\label{nonsplitting}

In this section, we give an example of a 
toric Fano manifold of index $2$ 
which admits no projective space bundle structure. 

Let $X=X_{\Sigma}$ be a smooth complete toric $d$-fold. 
For any $x\in\G(\Sigma)$ and $p\in{\mathbb Z}_{\geq 2}$, 
we construct a new toric manifold $\mathcal{H}_{(x,p)}(X)$ as follows. 

Put $\overline{N}:=N\oplus{\mathbb Z}^{p-1}$ and let 
$\{e_1,\ldots,e_d,e_{d+1},\ldots,e_{d+p-1}\}$ be the 
standard basis for $\overline{N}$. 
Put $z_1:=e_{d+1},\ldots,z_{p-1}:=e_{d+p-1}$ and 
$z_p:=x-(z_1+\cdots+z_{p-1})$. 
We define a fan $\overline{\Sigma}$ in $\overline{N}$ 
as follows: 
The maximal cones of $\overline{\Sigma}$ are 
$\sigma+{\mathbb R}_{\geq 0}z_{i_1}+\cdots+
{\mathbb R}_{\geq 0}z_{i_{p-1}}$, 
where $\sigma$ is any maximal cone in $\Sigma$ and 
$1\leq i_1<\cdots<i_{p-1}\leq p$. Namely, 
$$X_{\overline{\Sigma}}=
\mathbb{P}_X\left(
\mathcal{O}^{\oplus p-1}_X\oplus\mathcal{O}_X(D_x)
\right),$$
where $D_x$ is the toric prime divisor corresponding to $x$. 
Then, obviously, we have an extremal primitive relation 
$z_1+\cdots+z_p=x$ of $X_{\overline{\Sigma}}$. 
So, we obtain a smooth complete toric $(d+p-1)$-fold 
$\mathcal{H}_{(x,p)}(X)$ 
by the corresponding blow-down. It is obvious that 
the Picard number of $\mathcal{H}_{(x,p)}(X)$ is same as $X$. Moreover, 
$\mathcal{H}_{(x,p)}(X)$ has the 
following property:

\begin{prop}
The primitive collections of $\mathcal{H}_{(x,p)}(X)$ are
\begin{enumerate}
\item $P\in\PC(\Sigma)$, where $x\not\in P$, and
\item $(P\setminus\{x\})\cup\{z_1,\ldots,z_p\}$, where $P\in\PC(\Sigma)$ 
and $x\in P$.
\end{enumerate}
Moreover, if $x+x_1+\cdots+x_m=b_1y_1+\cdots+b_ny_n$ is a 
primitive relation of $X$, then 
$z_1+\cdots+z_p+x_1+\cdots+x_m=b_1y_1+\cdots+b_ny_n$ is a 
primitive relation of $\mathcal{H}_{(x,p)}(X)$, 
while if $x_1+\cdots+x_m=bx+b_1y_1+\cdots+b_ny_n$ is a 
primitive relation of $X$, then 
$x_1+\cdots+x_m=bz_1+\cdots+bz_p+b_1y_1+\cdots+b_ny_n$ is a 
primitive relation of $\mathcal{H}_{(x,p)}(X)$.

\end{prop}

\begin{proof}
The primitive collections of $\overline\Sigma$ are 
the primitive collections of $\Sigma$ and $\{z_1,\ldots,z_p\}$. 
Then, we can calculate the primitive collections of 
$\mathcal{H}_{(x,p)}(X)$ easily (see Corollary 4.9 in \cite{sato1}).
\end{proof}

Now, we can describe an example of a 
toric Fano manifold of index $2$ 
which admits no projective space bundle structure. 

\begin{ex}\label{hanachangum}
Let $X=X_{\Sigma}$ be the del Pezzo surface of degree $7$. 
The primitive relations of $\Sigma$ are 
$x_1+x_3=x_2$, $x_1+x_4=0$, $x_2+x_4=x_3$ 
$x_2+x_5=x_1$ and $x_3+x_5=0$, where $\G(\Sigma)=\{x_1,x_2,x_3,x_4,x_5\}$. 
Put $$Y=Y_{\widetilde{\Sigma}}
:=\mathcal{H}_{(x_1,2)}\left(\mathcal{H}_{(x_2,2)}
\left(\mathcal{H}_{(x_3,2)}\left(
\mathcal{H}_{(x_4,2)}\left(\mathcal{H}_{(x_5,2)}\left(X
\right)\right)\right)\right)\right).$$
Then, the primitive relations of $\widetilde{\Sigma}$ are 
$$x_1+x'_1+x_3+x'_3=x_2+x'_2,\ 
x_1+x'_1+x_4+x'_4=0,$$
$$x_2+x'_2+x_4+x'_4=x_3+x'_3,\ 
x_2+x'_2+x_5+x'_5=x_1+x'_1\mbox{ and}$$
$$x_3+x'_3+x_5+x'_5=0,$$
where $\G(\widetilde{\Sigma})=
\{x_1,x_2,x_3,x_4,x_5,x'_1,x'_2,x'_3,x'_4,x'_5\}$. We remark that 
$Y$ is a smooth toric Fano $7$-fold of index $2$, the Picard number 
of $Y$ is $3$ and $Y$ has no projective space bundle structure.
\end{ex}

\begin{rem}
Similarly as in Example \ref{hanachangum}, 
for any $p\in{\mathbb Z}_{\geq 2}$, 
we can construct a toric Fano manifold of index $p$ which 
has no projective space bundle structure.
\end{rem}

\end{document}